\documentclass[12pt]{article}

\usepackage{graphicx}
\usepackage{latexsym,amssymb}
\usepackage{amsthm}
\usepackage{indentfirst}
\usepackage{amsmath}
\usepackage{color}
\usepackage{xcolor}
\usepackage{fourier}
\usepackage[all]{xy}
\usepackage[colorlinks=true]{hyperref}

\newtheorem{theoremalph}{Theorem}

\newtheorem*{Main Theorem}{Main Theorem}

\newtheorem{Theorem}{Theorem}[section]
\newtheorem*{Theorem A}{Theorem A}
\newtheorem*{Theorem A'}{Theorem A'}
\newtheorem*{Theorem B'}{Theorem B'}

\newtheorem{Proposition}[Theorem]{Proposition}
\newtheorem{Lemma}[Theorem]{Lemma}

\newtheorem{Question}{Question}

\newtheorem{Remark}[Theorem]{Remark}

\newtheorem{Corollary}[Theorem]{Corollary}

\newtheorem{Claim-numbered}[Theorem]{Claim}

 \def\NN{{\mathbb N}} 

 \def\RR{{\mathbb R}} 
\def\TT{{\mathbb T}}

 \def\ZZ{{\mathbb Z}}

\def\La{\Lambda}

    \def\cS{{\cal S}}
    
    \def\cU{{\cal U}}
   \def\cP{{\cal P}} 
    
\def\cF{{\cal F}}

\newcommand{\diff}{{\operatorname{Diff}}}

\def\diff{\operatorname{Diff}}

\def\diam{\operatorname{Diam}}

\def\supp{\operatorname{Supp}}
\def\ud{\operatorname{d}}
\def\e{{\varepsilon}}

\begin{document}

\title{Disintegrations of non-hyperbolic ergodic measures along the center foliation of DA maps
}

\author{Ali Tahzibi and  Jinhua Zhang\footnote{A. Tahzibi is supported by FAPESP (2017-06463-3) and CNPq. J.Zhang is partially supported by National Key R$\&$D Program of China (2021YFA1001900) and NSFC 12001027.}}


\maketitle

\begin{abstract}
 
We show that each non-hyperbolic ergodic measure of a partially hyperbolic diffeomorphism on $\TT^3$ which is homotopic to Anosov admits a full measure subset which intersects each center leaf in at most two points. 

\hspace{-1cm}\mbox
\smallskip

\noindent{\bf Mathematics Subject Classification (2010).} 37D30, 37A05, 37C05. 
\\
{\bf Keywords.}   Partial hyperbolicity, non-hyperbolic ergodic measure, Anosov diffeomorphism.
\end{abstract}


\section{Introduction}
\subsection{Background}

One attempt of research in smooth ergodic theory is to verify up to which degree the known ergodic features of hyperbolic dynamics appear in beyond uniformly hyperbolic dynamics. There are several forms of relaxing uniform hyperbolicity condition.
Partial hyperbolicity  had attracted many attentions. Classification of partially hyperbolic diffeomorphisms and proving ergodicity of volume measure are two main research topics. The following statement of Pugh and Shub: "a little hyperbolicity goes a long way toward ergodicity", has been a driving force for many ergodicity results in the partially hyperbolic setting.

In this paper we would like to address hyperbolicity in terms of global measure theoretical properties. Observe that a non-trivial center bundle is the obstruction to uniform hyperbolicity of genuinely partially hyperbolic diffeomorphisms. In our quest for ``a little hyperbolicity" we calculate the Lyapunov exponent along center bundle. If the center Lyapunov exponent is non-zero almost every where we are in the setting of so-called nonuniform hyperbolicity and several methods including Pesin theory are available. So, vanishing center Lyapunov exponent is one of  the main obstructions to understand ergodic properties of partially hyperbolic diffeomorphisms. 

A  naive way of getting rid of center bundle is to study the quotient space of center foliation (whenever exists). However, in general this quotient space is not even a Hausdorff topological space. By the way, one may consider disintegration of  volume along the center foliation and it turns out that in some cases the Lebesgue measure disintegrates into Dirac measure (or finitely many Dirac) along the global center leaves! 
Shub and Wilkinson~\cite{SW} considered $A \times Id : \mathbb{T}^3 \rightarrow \mathbb{T}^3$ where $A$ is a linear Anosov diffeomorphism of $\mathbb{T}^2$ and after perturbation found open sets of volume preserving partially hyperbolic dynamics with non-vanishing center Lyapunov exponent.  Using compactness of center foliation and non vanishing  center Lyapunov exponent, they proved that the center foliation is non-absolutely continuous, and in fact the  volume disintegrates into atomic measures along center foliation~\cite{RW}. See also Pesin and Hirayama \cite{PH}  who proved atomic disintegration in higher dimensional center case where they assumed sum of the center exponents is non-zero.

 Here we would like to address this phenomenon of Dirac disintegration in the setting of partially hyperbolic diffeomorphisms on  $3$-manifolds. Sometimes it is  called \emph{virtual hyperbolicity} (This term was coined by E. Lindenstrauss and K. Scmidt in \cite{LiSch}): the existence of a full volume subset which intersects every center leaf  in at most finitely many points (or orbits). This ``little (virtual) hyperbolicity" may  be useful to prove ergodic properties.  For instance,  the  Dirac disintegration of Lebesgue along center foliation of derived from Anosov diffeomorphisms plays an important role to prove Bernoulli property in \cite{PTV1}, and  the examples of  Dirac disintegration of volume along center foliation of derived from Anosov diffeomorphisms enabled the authors of \cite{PTV} to construct examples of minimal yet measurable foliations .


In this paper,  we will consider non-hyperbolic ergodic measures and study their disintegration along the center foliation. 

There are many deep results in the active area:  the classification of partially hyperbolic diffeomorphisms of $3$-manifolds. In the spirit of these classification results we highlight following three  classes:
\begin{itemize}
		\item fibered maps,
	\item perturbation of time-one map of Anosov flows,
	\item  derived from Anosov diffeomorphisms.
\end{itemize}

By the results in the works of Avila, Viana and Wilkinson (\cite{AVW1}, \cite{AVW2}), in the first two categories there is a ``Dirac versus Lebesgue" dichotomy. 
R. Var\~{a}o \cite{V} gave an example to show that in the derived from Anosov case, the disintegration along center foliation may be neither atomic nor Lebesgue. 
It is important to note that  in the first two categories, if the center Lyapunov exponent is non-zero then we are in the atomic disintegration case of dichotomy. However, there are examples with vanishing center Lyapunov exponent but Dirac disintegration along the center foliation (See examples in \ref{s.time-one-map} and \ref{circlecenter}).  

Our main result is in the  derived from Anosov setting where we prove that vanishing center exponent always yields atomic disintegration. In fact,  this result holds for any invariant ergodic probability measure. We emphasize that in particular this shows a main difference between our result and previous disintegration techniques where smoothness of measure is crucial.
 
Besides its intrinsic interest we use this result in some class of derived from Anosov examples to show that the so called invariance principle \cite{Fur,L,AV} does not hold in this context.

In this way we shed more light on the atomic disintegration along the center foliation of all three partially hyperbolic categories above (See \ref{allthree}).

\subsection{Statement of the results}

Recall that  a diffeomorphism $f\in\diff^1(M)$ is \emph{partially hyperbolic} if there exists  a $Df$-invariant  splitting $TM=E^s\oplus E^c\oplus E^u$ and $N\in\mathbb{N}$ such that
$$\|Df^N|_{E^s(x)}\|<\min\{1,m(Df^N|_{E^c(x)}) \}\leq \max\{1,\|Df^N|_{E^c(x)}\| \}<m(Df^N|_{E^u(x)}).$$
A diffeomorphism $f\in\diff^1(\TT^3)$ induces an isomorphism $f_*: H_1(\TT^3,\ZZ)\to H_1(\TT^3,\ZZ)$ which can be  considered as a matrix in $GL(3,\ZZ)$, and one says that $f$ is \emph{derived from Anosov} (or \emph{homotopic to Anosov}) if $f_*$ induces an Anosov automorphism on $\TT^3.$

The existence of non-hyperbolic ergodic measures for robustly non-hyperbolic systems has been extensively investigated, see for instance~\cite{GIKN, KN, BBD1, BZ}. 
 Now, we investigate the disintegration of non-hyperbolic ergodic measures for derived from Anosov systems on $\TT^3$ which is one of the classical partially hyperbolic models and whose study originated from Ma\~n\'e~\cite{M}. 

\begin{theoremalph}~\label{thm.main-thm}
	Let $f\in\diff^1(\TT^3)$ be a partially hyperbolic diffeomorphism which is derived from Anosov. For any non-hyperbolic ergodic measure $\nu$, there exists a $\nu$-full measure set $\Lambda_\nu$  which intersects each center leaf in at most two points. 
\end{theoremalph}

\begin{Remark}
	\begin{enumerate}
		\item 	It has been shown in \cite{VY} that if the linear Anosov diffeomorphism has 2-dimensional unstable bundle, then each ergodic measure with negative center Lyapunov exponent admits a full measure set which intersects each center leaf in at most one point.
		\item Note that the center foliation of $f$ is orientable. If $f$ preserves the orientation of the center foliation,  the $\nu$-full measure set $\Lambda_\nu$   intersects each center leaf in at most one point. If the $\nu$-full measure set $\Lambda_\nu$   intersects almost every center leaf in two points, then the conditional measures of $\nu$ along the center foliation are equi-distributed on these two points. 
		\item In terms of  vanishing of transverse entropy (defined in \cite{LY2}), F. Ledrappier and J. Xie \cite{LX} showed that for $C^2$-diffeomorphisms, if the transverse entropy of an ergodic measure vanishes, then the conditional measure along a typical  unstable manifold is actually  carried  by a single strong unstable manifold. Thus the disintegrations along the weak unstable manifolds are atomic.  While we work for the case where the center Lyapunov exponent of a $C^1$-partially hyperbolic diffeomorphism vanishes, so our result is not implied by theirs. 
	\end{enumerate}
	
\end{Remark}
In \cite{GS}, the authors show that $C^{1+\alpha}$ volume preserving  partially hyperbolic diffeomorphisms $f$  which are derived from   Anosov on $\mathbb{T}^3$ are ergodic, and notice that  the Lebesgue measure is ergodic for $f^k$ ($k\in\ZZ\setminus\{0\}$).  Thus, one has the following corollary. 
\begin{Corollary}~\label{coro.main-thm}
	Let $f\in\diff^{1+\alpha}_m(\TT^3)$ be a partially hyperbolic  diffeomorphism which is derived from Anosov. If the center Lyapunov exponent of $f$ vanishes, then there exists a Lebesgue-full measure set  which intersects each center leaf in at most one point. 
\end{Corollary}

\subsection{Atomic disintegration and vanishing center exponent} \label{allthree}

In this section we  first mention some byproducts of our main result for derived from Anosov diffeomorphisms. Then we give some examples of vanishing center Lyapunov  exponent and atomic disintegration in the case of partially hyperbolic diffeomorphisms which are close to time one map of geodesic flow on negatively curved surface or  have compact center foliations.
 \subsubsection{Derived from Anosov diffeomorphisms}  

As we show in Theorem~\ref{thm.main-thm},  the disintegration of any ergodic measure along center foliation is Dirac if the center Lyapunov exponent vanishes.

Consider a family of examples in \cite{PT} where the authors begin with an Anosov automorphism and make Baraviera-Bonatti \cite{BB} type perturbation carefully to obtain derived from Anosov diffeomorphisms with some nice properties. 
More precisely, they begin with a linear Anosov automorphism (chosen from a suitable family of automorphisms) $A : \mathbb{T}^3 \rightarrow \mathbb{T}^3$ with eigenvalues $\lambda^s < \lambda^c < 1 < \lambda^u$ where $\lambda^c$ is very close to one,  and after   $C^1$-small perturbation they find a partially hyperbolic volume preserving diffeomorphism $f$  with $\lambda^c (f, m) > 0.$ It is proved in \cite{PT} that the disintegration of Lebesgue measure along center foliation of $f$ is atomic.  As the perturbation is done along a path, in fact the authors also obtain  a partially hyperbolic diffeomorphism $f_0$ homotopic to $A$ with vanishing center Lyapunov exponent.  
By Corollary \ref{coro.main-thm}, the disintegration of $m$ along center foliation of $f_0$ is given by Dirac measures. 

The perturbation in the above mentioned example is done locally and inside center-unstable bundle of Anosov automorphism. 
An immediate corollary is that the center-unstable foliation of $f$ and $A$ are the same and consequently all the leaves are planes and  it is an absolutely continuous foliation (in fact smooth).

\begin{Question}
	Is there any example of volume preserving derived from Anosov diffeomorphism $f$ such that the sign of its center Lyapunov exponent agrees with the sign of the center Lyapunov exponent of its linearization and the disintegration of Lebesgue measure along center foliation is Dirac?
\end{Question}

We will use this example and our theorem to show that the invariance principle in the non-compact center leaves does not hold here. (talking about invariance principle)

Crovisier and Poletti \cite{CP} obtained  the following invariance principle result. They consider partially hyperbolic and dynamically coherent diffeomorphisms $f$ which act quasi-isometrically in the center: there exist  $K \geq 1$ and $q > 0$ such that for every $x, y \in M$ with $y \in \mathcal{F}^c(x)$ and every $n \in \mathbb{Z},$ one has 
$$
K^{-1} d^c(x, y) - q \leq d^c(f^n(x), f^n(y)) \leq  K d^c(x, y) + q,
$$
where $d^c(\cdot,\cdot)$ denotes the distance along the center leaves.

\begin{Theorem}[\cite{CP}]
	Let $f\in\diff^1(M)$ be a partially hyperbolic and dynamically coherent diffeomorphism,  and let $m$ be an ergodic measure. Assume, in addition, that 
	\begin{itemize}
		\item  $f$ acts quasi-isometrically in the center;
		\item  $m$ has local $cu \times s$-product structure and $\lambda^c (f, m) = 0$.
	\end{itemize}
Then there exists a family of local center measures $\{m^c_x\}_{x \in \supp(m)}$ which are continuous, $f$-invariant, $s$-invariant, $u$-invariant and extend the center disintegration of $m$ to the whole support. 
\end{Theorem}

As we mentioned above, we can construct a derived from Anosov diffeomorphism $f$ such that center-unstable foliation is smooth and  this implies that volume measure $m$ locally has product structure. 

Now we claim that the conclusions of the above theorem can not hold for such an example. Indeed, by theorem \ref{thm.main-thm} the center disintegration is atomic. So take any such atom $x \in \mathbb{T}^3$ and consider a small neighborhood of $x$. By accessibility\footnote{Recall that a partially hyperbolic diffeomorphism is \emph{accessible} if any two points can be joined by a path which is a concatenations of finitely many curves lying either in a strong stable leaf or in a  strong unstable leaf.} (here we have  local non joint integrability) we can take an $su-$path beginning from $x$ and ending at some $y \neq x$ in the center plaque of $x$. As  $m$ is fully supported, if there exists a continuous disintegration of $m$ which is both $s$ and $u$-invariant then we would conclude that $y$ is also an atom of the disintegration of $m$ along the center plaque of $x$ which is a contradiction.
\subsubsection{Perturbation of time-one map}~\label{s.time-one-map}
Perturbations of time$-1$ map of geodesic flow on the unit  tangent bundle of negatively curved surfaces are classical example of partially hyperbolic diffeomorphisms isotopic to identity.  In this case a deep result of Avila-Viana-Wilkinson proves a dichotomy (Theorem \ref{AVW1}) : either the center foliation is virtually hyperbolic or the dynamics is embedded into a flow. 
We point out that this dichotomy is not a dichotomy in the level of center Lyapunov exponent. More precisely, if the center Lyapunov exponent is non-zero the authors show that there is a full measure subset which intersects almost every leaf in exactly $k-$orbits, but if the center Lyapunov exponent vanishes, the center foliation may not  be absolutely continuous. If the center foliation is absolutely continuous,  they prove a rigidity result (See item (2) of Theorem \ref{AVW1}). However, it is not difficult to see there are examples of vanishing center exponent and Dirac center disintegration. 

\begin{Theorem}[Main Theorem 1 in \cite{AVW1}]( Dirac - Lebesgue dichotomy)  \label{AVW1}
	Let $\phi^t : T^1 S \rightarrow T^1S$
	be the geodesic flow for a  negatively curved closed surface $S$ and let $m$ be the $\phi^t$-invariant Liouville probability
	measure.
	
	Then there is a  $C^1$-open neighborhood $\mathcal{U}\subset \diff^1_m(T^1S)$ of $\phi^1$   such that for any smooth diffeomorphism $f \in \cU$, one has 
	\begin{enumerate}
		\item either there exist $k \geq 1$ and a full $m$-measure set $Z \subset T^1S$ that intersects
		every center leaf in exactly $k$ orbits of f,
		\item or $f$ is the time-one map of an $m$-preserving $C^{\infty}$ flow.
	\end{enumerate}
		In the first case,  $m$ has atomic disintegration, and in the second case, it has Lebesgue
disintegration along the center foliation of $f$.
\end{Theorem}

The study of Lyapunov exponents and the so called Invariance Principle is crucial in the proof of the above theorem. By the way the above dichotomy does not coincide with hyperbolic versus non-hyperbolic behavior in the center direction. If  $m$ is non-uniformly hyperbolic (the center Lyapunov exponent is non-zero),  then   we are in the first item of dichotomy. However in the presence of zero Lyapunov exponent we may still be in the first case. 

We construct a family $\{g_t\}_{t \in [-1, 1]}\subset\diff^1_{m}(T^1S)$ of diffeomorphisms in a $C^1$ small neighborhood of $\phi^1$  such that $g_0 = \phi^1$ and $g_1$ (resp. $g_{-1}$) has positive (resp. negative) center Lyapunov exponent.

Indeed, by Baraviera-Bonatti \cite{BB} local perturbation method  consider two families of  volume preserving diffeomorphisms $(\xi_t)_{t \in [0,1]}$ and $(\eta_t)_{t \in [-1,0]}$ such that 
\begin{itemize}
	\item $\eta_0 = Id$ and  $\lambda^c(g_0 \circ \eta_1) > 0$,
	\item $\xi_0 = Id$ and  $\lambda^c(g_0 \circ \xi_{-1}) < 0$,
	\item  the supports of all $\xi_t$ and $\eta_t$ are in some small disjoint balls  $B_1$ and  $B_2$ respectively.  
\end{itemize}
Now define 
$$g_t = \begin{cases}
	g_0 \circ \eta_t  & t \in [0, 1]\\
	g_0 \circ \xi_t &  t \in [-1, 0].
\end{cases}$$
Let $\mathcal{C}$ be  a periodic orbit of $\phi^t$  which is  disjoint from $B_1$ and $B_2$. Then $\mathcal{C}$ is normally hyperbolic for $g_0$, and thus  admits continuation under small perturbations. By the conservative version of Franks lemma (see Proposition 7.4 in \cite{BDP}), there exists $g$ arbitrarily $C^1$ close to $g_0$ such that $g$ admits a sink on the compact center leaf $\mathcal{C}_g$ (continuation of $\mathcal{C}$ for $g$) and $g$ coincides with $g_0$ in the complement of a small neighborhood of $\mathcal{C}$.
Finally we define $f_t = g\circ g_0^{-1} \circ g_t $. As the center bundle is one dimensional,  the center Lyapunov exponent varies continuously. Since $\lambda^c (f_1) > 0 $ and $ \lambda^c(f_{-1}) < 0$,   by mean value theorem there exists $t_{*}$ such that $\lambda^c(f_{t_{*}}) = 0$. Observe that $f_{t_{*}}$ restricted to $\mathcal{C}_g$ admits a sink and thus it can not be embedded into a flow as in Theorem~\ref{AVW1}.

\subsubsection{Fibered partially hyperbolic diffeomorphisms (Circle center leaves)} \label{circlecenter} 
An important class of partially hyperbolic diffeomorphisms are the so called fibered maps and they include skew products. A diffeomorphism $f: M \rightarrow M$ is fibered partially hyperbolic if $M$ admits an $f-$invariant structure $\pi: M \rightarrow B$ of continuous fiber bundle with $C^1$ fibers where the fibers are tangent to the center bundle of $f$ (See \cite{AVW2} for more details.) Here we assume the fibers are homeomorphic to circle. 
 
 In the seminal work~\cite{SW}, Shub and Wilkinson show that one can perturb  a volume preserving skew-product system such that the new system have   non-absolutely continuous center foliation and the center Lyapunov exponent is non-vanishing;     in fact,  the disintegration of the volume along the center foliation is atomic (see \cite{RW}). 

Katok's example shows that zero center Lyapunov exponent and Dirac disintegration can occur simultaneously. However this example is not accessible.
We recall Katok's example here. Let $\{f_t\}_{t \in \mathbb{S}^1}$ be a smooth family of volume preserving  Anosov diffeomorphisms of $\mathbb{T}^2$ such that for any $s \neq t$ one has that $f_s$ and $f_t$ are topologically but not smoothly conjugate. By a result of de la Llave \cite{dL} the conjugacy between $f_s$ and $f_t$ can not be absolutely continuous. 
Now define $F : \mathbb{T}^2 \times \mathbb{S}^1 \rightarrow \mathbb{T}^2 \times \mathbb{S}^1$ by $F(x, t) = (f_t(x), t).$ It is not difficult to see that $\lambda^c (F, vol) = 0$ and the center holonomy between two transverse torus $\mathbb{T}^2 \times \{s\}$ and $\mathbb{T}^2 \times \{t\}$ is given by the conjugacy between $f_t$ and $f_s$ and hence not absolutely continuous. In fact considering the union of generic points of $f_t$ for $t \in \mathbb{S}^1$,  one obtains a full Lebesgue measure subset of $\mathbb{T}^3$ which intersects any center leaf in at most one point. 

Using a similar construction as in previous subsection,  one may find examples of  accessible volume preserving partially hyperbolic diffeomorphism on $M^{3}$ with vanishing center Lyapunov exponent and Dirac center disintegration. Take $M^3$ a nilmanifold on which   all the partially hyperbolic diffeomorphisms are accessible \cite{HHU}. Take a one parameter family of volume preserving partially hyperbolic diffeomorphisms on $M^3$ as in Section~\ref{s.time-one-map} 
and after a perturbation introduce a sink in a periodic center leaf of all diffeomorphisms in the path. There exists a parameter with vanishing center Lyapunov exponent and by rigidity result of \cite{AVW2} (Theorems A and B) the conditional measures along center foliation are atomic.


\medskip

\noindent{\bf Acknowledgments: } We would like to thank S. Crovisier, A. Wilkinson, S. Gan, M. Poletti and Y. Shi for helpful comments.

\section{Preliminaries}
In this section, we collect the notions and results used in this paper.
\subsection{Rohlin's Disintegration Theorem}
Given a compact metric space  $X$, a partition $\cP$ of $X$ is called \emph{measurable}, if there exists a sequence  of finite (Borel) measurable partitions  $\cP_1\prec\cP_2\prec\cdots\prec\cP_n\prec\cdots$ such that   $\cP=\vee_{n\in\NN} \cP_n$.

	Given a probability measure $\mu$ and a measurable partition $\cP$ on $X$, we denote by $\hat\mu$ the quotient measure  induced by $\mu$ on the Lebesgue space  $X/\cP$. 
\begin{Theorem}[\cite{Ro}]
	Given a probability measure $\mu$ and a measurable partition $\cP$ on $X$,  there exists a (essentially) unique family of probability measures $\{\mu_{P}\}_{P\in\cP}$ such that 
\begin{itemize}
	\item $\mu_P(P)=1$ for $\hat\mu$ a.e. $P\in\cP$;
	\item  for any measurable set $A\subset X$, the map $P\mapsto\mu_{P}(A)$ is measurable and 
	\[\mu(A)=\int\mu_P(A)\ud\hat\mu.
	\] 
\end{itemize}
	\end{Theorem}

The family of probability measure $\{\mu_{P}\}_{P\in\cP}$ is called the \emph{conditional measures or the disintegrations} of $\mu$ with respect to $\cP$.
\subsection{Partially hyperbolic diffeomorphisms homotopic to Anosov}
In this section, we recall some  properties of partially hyperbolic diffeomorphisms on $\TT^3$ which are homotopic to Anosov.

 \begin{Theorem}[\cite{F}]\label{thm.Franks}
 	Let $f$ be a homeomorphism on $\TT^3$, which is homotopic to a linear Anosov diffeomorphism $A$. Consider a lift $F$ of $f$ to $\TT^3$. Then there exists a unique continuous surjective map $\pi:\RR^3\to\RR^3$ such that 
 	\begin{itemize}
 		\item $\pi\circ F=A\circ \pi;$
 		\item $\pi(x+n)=\pi(x)+n$, for any $x\in\RR^3$ and $n\in\ZZ^3.$
 	\end{itemize}
 \end{Theorem}
\begin{Remark}
$\pi$ is homotopic to identity and 	$\pi$ induces a continuous surjective map on $\TT^3$, still denoted by $\pi$,  such that $\pi\circ f=A\circ \pi$.
\end{Remark}

\begin{Theorem}[\cite{BI, Ha,U,Po}]\label{thm.semi}
	Let $f$ be a partially hyperbolic diffeomorphism on $\TT^3$ which is homotopic to a linear  Anosov diffeomorphism $A$. Let $\pi$ be the semi-conjugacy between $f$ and $A$. Then one has the following properties:
	\begin{itemize}
		\item there exist unique foliations $\mathcal{F}^{cs}$, $\mathcal{F}^{cu}$ and $\cF^c$ tangent to the bundles $E^s\oplus E^c$, $E^c\oplus E^u$ and $E^c$ respectively;
		\item  the lift of $\cF^c$ to $\RR^3$ is quasi-isometric\footnote{A foliation$\cF$ with $C^1$-leaves on $\RR^3$ is quasi-isometric, if there exists $a,b>0$ such that $\ud_{\cF}(x,y)<a\ud(x,y)+b$ for any $x\in\RR^3$ and $y\in\cF(x)$.};
		\item $A$ has simple spectrum;
		\item the semi-conjugacy sends the center foliation of $f$ to the center foliation of $A$;
		\item  for an point $x\in\TT^3$, the set $\pi^{-1}(x)$ is a center closed segment (including a single point case).
	\end{itemize}
\end{Theorem}

\subsection{Entropy along an expanding foliation}
The entropy along unstable manifolds was defined in \cite{LY1, LY2} and is generalized for an expanding foliation in \cite{VY,Y,HHW}.

Let $f\in\diff^1(M)$.  A measurable partition $\xi$ is \emph{increasing} if  $f\xi\prec \xi$. 
A measurable partition $\xi$ is called  $\mu$-subordinate to a foliation $\cF$ on $M$ if 
  for $\mu$ a.e. $x$, the element $\xi(x)$ is contained in the leaf $\cF(x)$ and there exists $\delta_x>0$ such that the $\delta_x$ neighborhood $\cF_{\delta_x}(x)$ of $x$ (with respect to the leaf topology) is contained in $\xi(x)$. 
 
 A foliation $\cF$ on $M$ is called an \emph{expanding foliation} of $f\in\diff^1(M)$, if 
 \begin{itemize}
 	\item $\cF$ is $f$-invariant;
 	\item  each leaf of $\cF$ is $C^1$ and $f$ is uniformly expanding along the tangent bundle of    $\cF$.
 \end{itemize}
  
\begin{Lemma}[Proposition  3.1 in \cite{LS} and Lemma 3.2 in~\cite{Y}]~\label{l.increasing-partition}
	Let $f\in\diff^1(M)$,  $\mu$ be an invariant measure and $\cF$ be an expanding foliation $f$. Then there exists an increasing measurable partition which is $\mu$-subordinate to $\cF$. 
\end{Lemma}
Let $\mu$ be an invariant measure of a partially hyperbolic diffeomorphism $f$, then \emph{the unstable (metric) entropy} of $\mu$ is defined as 
$$h_\mu(f,\cF^u):= H_\mu(\xi|f\xi),$$
where $\xi$ is an increasing partition $\mu$-subordinate to the unstable foliation $\cF^u$. It has been shown in \cite{LY1} that the unstable entropy is independent of the choice of $\xi$.

The following results show that the disintegration of an ergodic measure along an expanding foliation is closely related to its metric entropy along this foliation. 
\begin{Proposition}[Proposition 2.5 in \cite{VY}] \label{p.zero-entropy}
	Let $f\in\diff^1(M)$, $\nu$ be an ergodic measure and  $\cF$ be an expanding foliation. Then the followings are equivalent:
	\begin{itemize}
		\item $h_\nu(f,\cF)=0.$
		\item  there exists a $\nu$-full measure subset intersecting  each center leaf in at most one point.
	\end{itemize}
\end{Proposition}
\begin{Proposition}[Proposition 2.7 in \cite{VY}] \label{p.positive-entropy}
	Let $f\in\diff^1(M)$, $\nu$ be an ergodic measure and  $\cF$ be an expanding foliation. Then the followings are equivalent:
	\begin{itemize}
		\item $h_\nu(f,\cF)>0.$
		\item  any   $\nu$-full measure subset intersects almost every center leaf in an uncountable set.   
	\end{itemize}
\end{Proposition}
\section{Entropy and disintegrations along the center foliation: Proof of Theorem~\ref{thm.main-thm}}
In this section, we first show that the projection of a non-hyperbolic ergodic measure has zero entropy along the center foliation, then we study the disintegration along the center foliation.  

\begin{Proposition}~\label{p.key-proposition}
	Let $f\in\diff^1(\TT^3)$ be a partially hyperbolic diffeomorphism homotopic to a linear Anosov diffeomorphism $A$. Let $\pi:\TT^3\to\TT^3$ be the semi-conjugacy between $f$ and $A$ given by Theorem~\ref{thm.Franks}. For any non-hyperbolic ergodic measure $\nu$, one has $h_{\hat\nu}(A, W^c)=0$ where $\hat\nu:=\pi_{*}\nu$.
\end{Proposition}
\proof  
Without loss of generality, one can assume that $A$ is expanding along the center. Let us denote by $\cF^c$ and $W^c$ the center foliation of $f$ and $A$ respectively. Let $\nu$ be an ergodic measure with vanishing center Lyapunov exponent, then $\hat\nu:=\pi_{*}\nu$ is an ergodic measure of $A$. 
The proof proceeds by contradiction. Assume, on the contrary, that $h_{\hat\nu}(A, W^c)>0.$ 
\begin{Claim-numbered}	~\label{c.isomorphism}
	 $(f,\nu)$ is isomorphic to $(A,\hat\nu)$ via $\pi$. 
	\end{Claim-numbered}
\proof[Proof of Claim~\ref{c.isomorphism}]
By Theorem~\ref{thm.semi}, the set $\pi^{-1}(x)$ is a center segment (could be trivial) for any $x\in\TT^3$.
Consider the set $$\cS=\big\{x\in\TT^3| \textrm{$\pi^{-1}(x)$ is a non-trivial center segment} \big\}.$$ Notice that $f^{-1}(\pi^{-1}(x))=\pi^{-1}(A^{-1}(x))$ for any $x\in\TT^3$. Thus, the set $\cS$ is $A$-invariant. 
By the ergodicity of $\hat\nu$, one has $\hat\nu(\cS)=0$ or $1$. As $h_{\hat\nu}(A, W^c)>0$, by Proposition~\ref{p.positive-entropy}, there exists a $\hat\nu$-full measure subset $K$ such that  $K\cap W^c(x)$ is an uncountable set for $\hat\nu$ a.e. $x$. Since on each   leaf of $W^c$ there are at most countable many points whose pre-images under $\pi$  are non-trivial center segments, thus $\hat\nu(\cS)=0$ proving that $\pi$ is an isomorphism between $(f,\nu)$ and $(A,\hat\nu)$.
\endproof

Now, consider a measurable partition $\hat\xi^{c}$  which is increasing partition  and $\hat\nu$-subordinate to the foliation $W^c$ with $\diam(\hat\xi^{c})<1$. Let $\{\hat{\nu}_x^c\}$ be the conditional measures of $\hat\nu$ with respect to the measurable partition $\hat\xi^c.$ By Shannon-McMillan-Breiman theorem for unstable entropy \cite[Lemma 9.3.1]{LY2},  one has 
$$\lim_{n\rightarrow\infty}-\frac{\log \hat\nu_x^{c}(\bigvee_{i=0}^{n-1}A^{-i}\hat\xi^{c}(x))}{n}=h_{\hat\nu}(A, W^c) \textrm{,\:\:  for $\hat\nu$ a.e. $x\in\TT^3$}.$$
As $\pi:\TT^3\to\TT^3$ is continuous, the partition $\xi^{c}:=\pi^{-1}(\hat\xi^{c})$  is measurable. 
Since  $\pi$ is homotopic to identity and the center foliation is quasi-isometric, there exists $\eta_0>0$ such that $\diam(\xi)<\eta_0$. 
By Theorem~\ref{thm.semi},  $\pi$ sends the center leaves of $f$ to the center leaves of $A$ and $\pi^{-1}(x)$ is a center segment for any $x\in\TT^3$, and thus the partition $\xi^c$ is $\nu$-subordinate to the center foliation $W^c$. Notice that $\xi^c=\pi^{-1}(\hat\xi^c)\prec\pi^{-1}A^{-1}(\hat\xi^c)=f^{-1}\pi^{-1}(\xi)=f^{-1}(\xi^c)$. To summarize, the partition $\xi^c$ is an increasing measurable partition $\nu$-subordinate to $\cF^c$. We denote by $\{\nu_x^{c}\}$ the conditional measures of $\nu$ with respect to the measurable partition $ \xi^{c}$.

Since $\pi$ is an isomorphism between $(f,\nu)$ and $(A,\hat\nu)$, one has   $$h_{\hat\nu}(A,W^c)=H_\nu(f^{-1} (\xi^{c})|\xi^{c})=\lim_{n\rightarrow\infty}-\frac{\log\nu_x^{c}(\bigvee_{i=0}^{n-1}f^{-i}\xi^{c}(x))}{n} \textrm{,\:\:  for $\nu$ a.e. $x\in\TT^3$}.$$   
Now fix $\e\in(0,\min\{h_{\hat\nu}(A, W^c)/20,1/8\})$  small, and let us denote 
$$\Lambda_{N,\e}=\big\{x\in\TT^3|\; \nu_x^{c}(f^{-n}(\xi^{c}(f^n(x))))\leq e^{-n(h_{\hat\nu}(A, W^c)-\e)}, \textrm{ for any $n\geq N$}\big\},$$ 
$$K_{N,\e}=\big\{x\in\TT^3|\;  e^{-|n|\e}\leq \log\|Df^n|_{E^c(x)}\|\leq e^{|n|\e}, \textrm{ for any $|n|\geq N$}\big\}.$$
Then $\nu(\Lambda_{N,\e}\cap K_{N,\e})$ tends to $1$ when $N$ tends to infinity for any $\e>0.$

For $\delta>0$, we define the set  $$L_\delta=\big\{x\in\TT^3| \cF^c_\delta(x)\subset \xi^{c}(x)\big\}.$$
As $\xi^{c}$ is $\nu$-subordinate to the center foliation $\cF^c$, $\nu(L_\delta)$ tends to $1$ when $\delta$ tends to $0$.
Choose $\delta>0$ small such that  $\nu(L_\delta)>3/4$, and up to shrinking $\delta$, one can assume that 
 $$e^{-\e}<\|Df|_{E^c(x)}\|/\|Df|_{E^c(y)}\|<e^\e \textrm{\:\: for any $x,y\in\TT^3$ with $\ud(x,y)<\delta.$} $$
We consider the set 
$$R_{N,\e}=\big\{x\in \TT^3|\;\;\nu(L_\delta)-\e< \frac 1n\sum_{i=0}^{n-1}\chi_{_{L_\delta}}(f^i(x))<\nu(L_\delta)+\e,\textrm{ for any $n\geq N$}\big\}.$$
By Birkhoff ergodic theorem, $\nu(R_{N,\e})$ tends to $1$ when $N$ tends to infinity.
Let $\tau_{n,\e}=[8n\e]+1$ for $n\in\NN.$
\begin{Claim-numbered}~\label{cl.return}
For any $n>N$  and  any $x\in R_{N,\e}$, there exists $j\in\{n,\cdots, n+\tau_{n,\e}-1 \}$ such that $f^j(x)\in L_\delta.$
\end{Claim-numbered}
\proof[Proof of Claim~\ref{cl.return}]
By the definition of $R_{N,\e}$, for $n>N$, one has 
$$\nu(L_\delta)-\e< \frac 1n\sum_{i=0}^{n-1}\chi_{L_\delta}(f^i(x))<\nu(L_\delta)+\e,$$
$$\nu(L_\delta)-\e< \frac{1}{n+\tau_{n,\e}}\sum_{i=0}^{n+\tau_{n,\e}-1}\chi_{L_\delta}(f^i(x))<\nu(L_\delta)+\e.$$
Thus for $n>N$, by the choice of $\e$, one has 
$$\sum_{i=n}^{n+\tau_{n,\e}-1}\chi_{L_\delta}(f^i(x))>(n+\tau_{n,\e})(\nu(L_\delta)-\e)-n(\nu(L_\delta)+\e)>\tau_{n,\e}(\nu(L_\delta)-\e)-2n\e>\tau_{n,\e}/2-2n\e>0.$$   
This ends the proof of Claim~\ref{cl.return}.
\endproof	

\bigskip

 Fix $N$ large enough such that there exists a subset $\tilde\La\subset\Lambda_{N,\e}\cap K_{N,\e}\cap R_{N,\e}\cap L_\delta$ satisfying that 
\begin{itemize}
	\item $\nu(\tilde{\La})>1/2$;
	\item $\nu_x^{c}(\Lambda_{N,\e}\cap K_{N,\e}\cap R_{N,\e}\cap L_\delta)>1/2$ for any $x\in\tilde{\La}.$
\end{itemize}

Recall that $\xi^c\prec f^{-1}(\xi^c)$. For each $x\in\tilde\La$ and $n>N$,  the elements of the partition  $f^{-n}\xi^{c}|_{\xi^{c}(x)}$\vspace{1.5mm}
 on $\xi^c(x)$ which contain  points in $\Lambda_{N,\e}\cap K_{N,\e}\cap R_{N,\delta}\cap L_\delta$ have $\nu_x^c$-measure at most $e^{-n(h_{\hat\nu}(A, W^c)-\e)}$,\vspace{1.5mm} and therefore there are at least 	$l_n:=[e^{n(h_{\hat\nu}(A, W^c)-\e)}/2]$ elements of $f^{-n}\xi^{c}|_{\xi^{c}(x)}$ intersecting $\Lambda_{N,\e}\cap K_{N,\e}\cap R_{N,\delta}\cap L_\delta$.   Let $x_1,\cdots,x_{l_n}\in \Lambda_{N,\e}\cap K_{N,\e}\cap R_{N,\delta}\cap L_\delta\cap \xi^{c}(x)$ satisfy that $$f^{-n}\xi^{c}(x_i)\cap f^{-n}\xi^{c}(x_j)=\emptyset  \textrm{\:\:  for $i\neq j.$}$$ 
By the definition of $L_\delta$, one has  $\cF^c_\delta(x_i)\subset \xi^{c}(x)$. Now, we decompose the set
$\{x_1,\cdots, x_{l_n}\}$ according to their return time to $L_\delta$.
 For each $j\in \{0,\cdots, \tau_{n,\e}-1 \}$, let $Z_j=\big\{x_i| f^{j+n}(x_i)\in L_\delta\big\}$.
By Claim~\ref{cl.return}, it holds $\cup_{j=0}^{\tau_{n,\e}-1}Z_j=\{x_1,\cdots, x_{l_n}\}$, then there exists  $j_0 \in\{0,\cdots, \tau_{n,\e}-1 \}$  such that $\#Z_{j_0}=\max\#Z_j\geq \frac{l_n}{\tau_{n,\e}}$. By the definition of $L_\delta$, one has 
 $$\cF^c_\delta(f^{j_0+n}(x_i))\subset \xi^{c}(f^{j_0+n}(x_i)) \textrm{\:\: for $x_i\in Z_{j_0}$.}$$

Let us denote  $\gamma_n=C^{-N}\cdot e^{-3n\e}\cdot\delta$, where $C=\sup_{x\in\TT^3}\|Df|_{E^c(x)}\|>1$.
\begin{Claim-numbered}~\label{c.center-size}
For any  $n>N$ and any $x_i\in Z_{j_0}$, one has 	$\cF^c_{\gamma_{_{n+j_0}}}(x_i)\subset f^{-j_0-n}(\xi^{c}(f^{j_0+n}(x_i)))$.
	\end{Claim-numbered}
\proof[Proof of Claim~\ref{c.center-size}] 
We shall inductively show that $f^j(\cF^c_{\gamma_{n+j_0}}(x_i))\subset \cF^c_{\delta}(f^j(x_i))$ for $j\leq j_0+n.$

For any $j\leq N$, one has that $$f^j(\cF^c_{\gamma_{n+j_0}}(x_i))\subset \cF^c_{C^j\gamma_{n+j_0}}(f^j(x_i))  \textrm{\: and \:}  C^j\cdot\gamma_{n+j_0}\leq e^{-3(n+j_0)}\delta<\delta.$$
Assume that $f^j(\cF^c_{\gamma_{n+j_0}}(x_i))\subset \cF^c_{\delta}(f^j(x_i))$ already holds for some $N\leq k\leq j_0+n$ and any $j\leq k$.  Then for any point $y\in\cF^c_{\gamma_{n+j_0}}(x_i)$, one has  $$\ud^c(f^{k+1}(y),f^{k+1}(x_i))=\|Df^{k+1}|_{E^c(y_{k+1})}\|\cdot \ud^c(y,x_i)\leq \|Df^{k+1}|_{E^c(y_{k+1})}\|\cdot\gamma_{n+j_0},$$  where   $y_{k+1}\in \cF^c_{\gamma_{n+j_0}}(x_i)$ and $\ud^c(\cdot,\cdot)$ denotes the distance along the center leaf.
Since $f^{j}(y_{k+1})\in\cF^c_\delta(f^j(x_i))$ for $0\leq j\leq k$, by the choice of $\delta$, one has that 
$$\|Df^{k+1}|_{E^c(y_{k+1})}\|=\prod_{j=0}^k\|Df|_{E^c(f^j(y_{k+1}))}\|\leq e^{(k+1)\e}\prod_{j=0}^k\|Df|_{E^c(f^j(x_i))}\|\leq e^{2(k+1)\e}.$$
Thus $\ud^c(f^{k+1}(y),f^{k+1}(x_i))\leq e^{2(k+1)\e}\cdot \gamma_{n+j_0}<\delta,$ proving that $f^{k+1}(\cF^c_{\gamma_{n+j_0}}(x_i))\subset \cF^c_{\delta}(f^{k+1}(x_i)).$
\endproof
As $f\xi^{c}\prec\xi^{c}$, by Claim~\ref{c.center-size}, for any $x_i\in Z_{j_0}$, one has $$\cF^c_{\gamma_{_{n+j_0}}}(x_i)\subset f^{-j_0-n}(\xi^{c}(f^{j_0+n}(x_i)))\subset f^{-n}(\xi^{c}(f^n(x_i))).$$
 By the choice of $x_i$, the element $\xi^{c}(x)$ contains at least $\#Z_{j_0}$-pairwise disjoint center segments of length $\gamma_{_{n+j_0}}$, and thus 
 $$\eta_0\geq \diam(\xi^{c}(x))\geq \#Z_{j_0}\cdot\gamma_{_{n+j_0}}\geq l_n\cdot\gamma_{_{n+j_0}}/\tau_{n,\e}\geq [e^{n(h_{\hat\nu}(A, W^c)-\e)}/2] C^{-N}\cdot e^{-3(n+[8n\e]+1)\e}\cdot\delta/([8n\e]+1).$$
 By the choice of $\e$, the right side tends to infinity when  $n$ tends to infinity,  which gives the contradiction. 
 \endproof

Now, we are ready to give the proof of our main result.
\proof[Proof of Theorem~\ref{thm.main-thm}]
Let $\pi:\TT^3\to\TT^3$ be the semi-conjugacy given by Theorem~\ref{thm.Franks}.  By Theorem~\ref{thm.semi}, the set $\pi^{-1}(\pi(x))$   is center segment (could be trivial) for any $x\in\TT^3.$
Since $\pi:\TT^3\to\TT^3$ is homotopic to identity and the center foliation is quasi-isometric, there exists a constant $K>0$ such that for any $x\in\TT^3$, the length of $\pi^{-1}(\pi(x))$ is bounded   from above by  $K$. For any $x\in\TT^3$, let us denote $\Gamma_x=\pi^{-1}(\pi(x))\subset\cF^c(x)$, then $f(\Gamma_x)=\Gamma_{f(x)}$.

Let $\nu$ be a non-hyperbolic ergodic measure, then $\hat\nu:=\pi_*(\nu)$  is an ergodic measure of $A$. By Proposition~\ref{p.key-proposition}, $h_{\hat\nu}(A, W^c)=0$. By Proposition~\ref{p.zero-entropy}, there exists a $\hat\nu$-full measure subset $\hat\La$ intersecting each center leaf in at most one point. Let $\La=\pi^{-1}(\hat\La)$ which is a $\nu$-full measure subset.  For any $x\in\La$, one has that $\La\cap \cF^c(x)= \Gamma_x$ which is a center segment (could be trivial). As $\pi$ is a continuous map, $\{\Gamma_x\}$ is a  measurable partition. 

 If there exists a $\nu$-positive measure set in which the center segment $\Gamma_x$ is trivial, then one can conclude by the ergodicity of $\nu$.

It remains the case where $\Gamma_x$ is  a non-trivial center segment for $\nu$ a.e. $x\in\TT^3$.  Let $\{ \nu_x^c\}$ be the conditional measures   of $\nu$ with respect to  the measurable partition $\{\Gamma_x\}$. As $f(\Gamma_x)=\Gamma_{f(x)}$ and $\nu$ is $f$-invariant, by the uniqueness of the conditional measures, one has $f_*\nu_x^c=\nu_{f(x)}^c$.  By the ergodicity of $\nu$, one has that
\begin{itemize}
	\item either $\nu_x^c$ has atoms for $\nu$. a.e. $x\in\TT^3$; or
	\item $\nu_x^c$ has no atoms  for $\nu$. a.e. $x\in\TT^3$.
\end{itemize}
We claim that the later case cannot happen.
\begin{Claim-numbered}~\label{c.admit-atom}
For $\nu$. a.e. $x\in\TT^3$, the conditional measure $\nu_x^c$ has atoms. 
\end{Claim-numbered}
\proof[Proof of Claim~\ref{c.admit-atom}]
Assume, on the contrary, that  $\nu_x^c$ has no atoms  for $\nu$. a.e. $x\in\TT^3$.	Notice that the center foliation $\cF^c$  is orientable and we fix an orientation. Up to replacing $f$ by $f^2$ and $\nu$ by its ergodic components under $f^2$; one can assume that $f$ preserves the orientation of $\cF^c$. We denote by  $\Gamma_x=[\alpha_x,\omega_x]^c$ such that  the orientation from $\alpha_x$ directed to $\omega_x$ gives the same orientation on $\Gamma_x$ as we fixed. 

For $\nu$. a.e. $x\in\TT^3$, one can consider the function $\phi^c_x:[\alpha_x,\omega_x]^c\to [0,1]$ defined by $\phi^c_x(z)=\nu^c_x([\alpha_x,z]^c)$. As the conditional measures have no atoms, the function $\phi^c_x$ is a continuous and non-decreasing function. Let $\beta_x\in[\alpha_x,\omega_x]^c$ be the point such that 
$\phi_x^c(\beta_x)=\frac{1}{2}$ and the length  of the center segment $[\alpha_x,\beta_x]^c$ is the smallest one. By the fact that $f(\Gamma_x)=\Gamma_{f(x)}$ and $f_*(\nu_x^c)=\nu_{f(x)}^c$, one has $f(\beta_x)=\beta_{f(x)}$. Now, we can consider the family of measures $\{2\cdot\nu^c_x|_{[\alpha_x,\beta_x]}\}$ which gives an invariant probability measure ${\mu}$ of $f$. By definition, the measure ${\mu}$ is absolutely continuous with respect to the $\nu$ and ${\mu}\neq\nu$ which contradicts with the ergodicity of $\nu$. 
\endproof

By Claim~\ref{c.admit-atom}, there exists an integer $k_0\in\NN$ such that   $$\nu\big(\big\{x\in\TT^3| \textrm{ there exists $z\in\Gamma_x$ with $\nu_x^c(z)\geq 1/k_0$} \big\}\big)>0.$$   
As $f_*\nu_x^c=\nu_{f(x)}^c$ and $f(\Gamma_x)=\Gamma_{f(x)}$, the set $\big\{x\in\TT^3| \textrm{ there exists $z\in\Gamma_x$ with $\nu_x^c(z)\geq 1/k_0$} \big\}$ is $f$-invariant, and thus by the ergodicity of $\nu$, one has 
 $$\nu\big(\big\{x\in\TT^3| \textrm{ there exists $z\in\Gamma_x$ with $\nu_x^c(z)\geq 1/k_0$} \big\}\big)=1.$$ 
Notice that the center foliation $\cF^c$ is orientable and now we fix an orientation of $\cF^c$. For each $\Gamma_x$, let $x^-\in\Gamma_x$ and $x^+\in\Gamma_x$ be the leftest and rightest points in $\Gamma_x$ with $\nu_x^c$-measure no less than $1/k_0$ (one could have $x^-=x^+$). If $f$ preserves the orientation of the center foliation, then $f(x^-)=(f(x))^-$ and $f(x^+)=(f(x))^+$ due to the fact that $f_*\nu_x^c=\nu_{f(x)}^c$. Consider the probability measure 
$\mu= \int\delta_{x^-} \ud\nu$, then $\mu$ is absolutely continuous to $\nu$ and the probability measure $\mu$ is $f$-invariant. By the ergodicity of $\nu$, one has that $\nu=\mu.$ Thus the disintegration of $\nu$ along the center foliation has exactly one point. 
If $f$ reverses the orientation of the center foliation,  then $f(x^-)=(f(x))^+$ and $f(x^+)=(f(x))^-$ due to the fact that $f_*\nu_x^c=\nu_{f(x)}^c$. 
Consider the probability measure 
$\mu=\frac{1}{2}\int(\delta_{x^-}+\delta_{x^+})\ud\nu$, then $\mu$ is absolutely and the probability measure $\mu$ is $f$-invariant. By the ergodicity of $\nu$, one has that $\nu=\mu.$ Thus the disintegration of $\nu$ along the center foliation has at most two points, and the conditional measures along the center foliation are equi-distributed on the atoms. 
\endproof

\begin{Remark}
	For a non-hyperbolic ergodic measure $\nu$, if its disintegration along the center foliation has two atoms, then $f$ must reverse the orientation of the center foliation and  $\nu=\frac{1}{2}(\nu_1+\nu_2)$ where $\nu_1,\nu_2$ are two different non-hyperbolic  $f^2$-ergodic measures, and the disintegrations of $\nu_1$ and $\nu_2$ along the center foliation have only one atom. By our proof,  $\pi_*(\nu_1)=\pi_*(\nu_2)$ and thus $\pi_*(\nu)$ is $A^2$-ergodic. 
\end{Remark}

 \bibliographystyle{plain}

\vspace{2mm}

 \begin{tabular}{l l l}
	\emph{\normalsize Ali Tahzibi}
	& \quad &
	\emph{\normalsize Jinhua Zhang}
	\medskip\\
	\small Instituto de Ci\^encias Matem\'aticas e de Computa\c{c}\~ao,
	&& \small School of Mathematical Sciences\\
	\small Universidade de S\~ao Paulo
	(USP)
	&& \small Beihang University\\
	\small  13566-590, S\~ao Carlos-SP, Brasil
	&& \small Beijing, 100191, P. R.  China\\
	\small \texttt{tahzibi@icmc.usp.br}
	&&\small \texttt{jinhua$\_$zhang@buaa.edu.cn}\\
	&&\small \texttt{zjh200889@gmail.com}
	
\end{tabular}

\end{document}